\documentclass[11pt]{article}
\usepackage{amsmath,amssymb,amsthm}
\usepackage[cmtip,arrow]{xy}
\usepackage{pb-diagram,pb-xy}
\usepackage{fancyhdr}

\newtheorem{theorem}{Theorem}[section]
\newtheorem{prop}[theorem]{Proposition}
\newtheorem{lem}[theorem]{Lemma}

\theoremstyle{definition}
\newtheorem{rema}[theorem]{Remark}

\newtheorem{sit}[theorem]{Situation}
\newtheorem{ex}[theorem]{Examples}
\newtheorem{defi}[theorem]{Definition}

\makeatletter
\let\c@equation=\c@theorem

\makeatother

\DeclareMathOperator{\lndet}{lndet}
\DeclareMathOperator{\tr}{tr}
\DeclareMathOperator{\id}{id}

\DeclareMathOperator{\pr}{pr}
\DeclareMathOperator{\im}{Im}

\DeclareMathOperator{\re}{Re}
\DeclareMathOperator{\Tr}{Tr}

\newcommand{\field}{\mathbb}
\newcommand{\Z}{\field{Z}}              
\newcommand{\N}{\field{N}}              
\newcommand{\R}{\field{R}}
\newcommand{\C}{\field{C}}
\newcommand{\Q}{\field{Q}}

\newcommand{\imp}{\Longrightarrow}

\renewcommand{\NG}{\mathcal{N}(G)}
\newcommand{\Np}{\mathcal{N}(\pi)}

\newcommand{\Ze}{\mathcal{Z}}
\newcommand{\NGI}{\mathcal{N}(G_i)}
\newcommand{\NU}{\mathcal{N}(U)}

\newcommand{\tru}{\tr^u_{\NG}}
\newcommand{\trg}{\tr^{\langle g\rangle}_{\NG}}
\newcommand{\trre}{\Tr^{\langle g\rangle,\re}_{\NG}}
\newcommand{\trrei}{\Tr^{\langle g\rangle,\re}_{\NGI}}
\newcommand{\trC}{\tr^\C_{\NG}}

\begin{document}

    \begin{center}
        {\huge Approximation of Center-Valued

        Betti-Numbers}
        \bigskip

        {Anselm Knebusch\footnote{Research supported by Marie Curie Research Training Network Non-Commutative Geometry MRTN-CT-2006-031962}
            at the Department of Mathematics

          Georg-August-Universit\"at G\"ottingen (Germany) and \\
          Katholieke Universiteit Leuven (Belgium)

          E-mail: knebusch@uni-math.gwdg.de}
          \bigskip
        \begin{abstract}\noindent
            A useful tool to calculate $L^2$-Betti-numbers is an approximation theorem which is
            proved in its original version by W. L\"uck. It shows that $L^2$-Betti-numbers $\beta^{(2)}_n(\widetilde{X})$ of the
            universal covering   $\widetilde{X}$ of a CW complex $X$\,,
            with residually finite fundamental group $\pi$\,, can be approximated by the Betti-numbers of the finite subcovers
            $\beta^{(2)}_n(\widetilde{X}/\pi_i)$\,. Since then, the the approximation theorem has been generalized in several steps
            and is now proven for a large class of groups containing e.g. all extensions of residually finite groups with
            amenable quotients, all residually amenable groups and free products of these.

            \noindent However, there is also a finer invariant than $L^2$-Betti-numbers: the so called universal or center-valued
            Betti-numbers $\beta^u$\,.
            They measure the dimension of the homology, using the center-valued trace $\tr^u_{\NG}$ of the finite von Neumann
            algebra $\NG$\,, instead of the usual $\C$-valued trace $\tr^\C_{\NG}$.
            In this paper we generalize the ordinary approximation theorem to an approximation theorem for universal Betti-numbers.
        \end{abstract}

    \end{center}

    \section{Introduction}
        For a finite CW complex $X$ with fundamental group  $\pi$\,, the $L^2$-homology of the universal covering
        $\widetilde{X}$ is given as the kernel of the combinatorial Laplacians $\Delta_*$ on
        $C_*^{(2)}(\widetilde{X})=C_*^{(cell)}(\widetilde{X})\otimes_{\Z \pi}\ell^{(2)}(\pi)$\,, which is after a choice of a
        cellular base isomorphic to a complex of finite direct sums of $\ell^2(\pi)$ on which the Laplacian
        $\Delta_p=(c_p\otimes \id)^*(c_p\otimes \id)+(c_{p-1}\otimes \id)(c_{p-1}\otimes \id)^*$ acts by left multiplication with
        a matrix over $\Z \pi \subset \Np$\,. Here, $\Np\subset\mathcal{B}(\ell^2(\pi))$ is the group von Neumann algebra of $\pi$:
        it is the von Neumann algebra generated by the left regular representation of $\pi$. $L^2$-Betti-numbers measure the
        dimension of the $L^2$-homology and can be defined as $\beta_p^2(X):=\dim_{\NG}^\C(\ker(\Delta_p))$\,.

        W. L\"uck shows in \cite{Luck1}  that the $L^2$-Betti-numbers $\beta^{(2)}_n(\widetilde{X})$ of the universal covering $\widetilde{X}$
        of a CW complex $X$\,, with residually finite fundamental group $\pi$\,, can be approximated by their finite
        dimensional analogons $\beta^{(2)}_n(\widetilde{X}/\pi_i)$\,.

        Using these ideas in a different context, J. Dodziuk and V. Mathai prove in \cite{D+M} a similar approximation result for amenable groups.
        In \cite{Schick}\,, T. Schick combines  both ideas and extends the result to a more general class $\mathcal{G}$, of groups, containing
        in particular amenable and residually finite groups.

        In an alternative formulation these proofs rely on showing that the kernel of a matrix $A\in M_d(\Z G)$ can be approximated via the kernels
        of the matrices $p_i(A)\in M_d(\Z G)$, where the $p_i$ are coming from some limit or extension process of $G$\,.
        Finally in \cite{Schick+} J. Dodziuk, P. Linnell, T. Schick and S. Yates, extend the coefficient ring $\Z G$ to $\overline{\Q} G$\,,
        especially to prove the Atiyah conjecture for matrices over $\overline{\Q} G$ and $G$ from a subclass of $\mathcal{G}$\,.
        In this paper, the approximation theorem will be generalized to an approximation theorem for the center-valued Betti-numbers
        $\beta_p^u(X):=\dim_{\NG}^u(\ker(\Delta_p))$\,.
        More precisely, we show that their Fourier coefficients (which are multiples of the so called delocalized Betti-numbers
        introduced by Lott in \cite{Lott}) can be approximated.

        \section{Notation and Preliminaries}

            We first need to introduce some notations. In the following, $G$ always denotes a discrete group,
            and we write $\mathcal{C}(G)$ for the set of elements $g\in G$ with finite conjugacy class $\langle g\rangle$.
            For abbreviation we denote by $g$ the elements $u_g\in\NG$\,.
            Using this notation the group ring is given by
            $\C G:=\{\sum_{g\in G}\lambda_gg\mid \text{only finitely many } \lambda_g\ne0\}\subset\NG$\,.
            The center of a von Neumann algebra $\mathfrak{A}$ is denoted by $\Ze(\mathfrak{A}):=\mathfrak{A}\cap \mathfrak{A}'$\,.
            The matrix ring $M_d(\NG)$ is defined as $M_d(\NG):=\NG\otimes_\C M_d(\C)$ and we let these
            operators act on $\ell^2(G)^d:=\ell^2(G)\otimes \C^d $\,.

            \begin{prop}
                The von Neumann algebra $\NG$ is a finite von Neumann algebra, hence it comes with some standard traces.
                \begin{enumerate}
                \item The universal or center-valued  trace
                \begin{align*}
                    \tru:\NG\longrightarrow \Ze(\NG)\,,
                \end{align*}
                which satisfies the following properties (see e.g. \cite{Kadison}):
                for all $a,b\in \NG$ and $c\in\Ze(\NG)$
                \begin{itemize}
                    \item $\tru(ab)=\tru(ba)$\,,
                    \item $\tru(c)=c$\,,
                    \item $\tru(a^*a)\ge 0$ \,\,\, and \,\,\,$\tru(a^*a)=0\imp a=0$\,,
                    \item $\tru(ca)=c\tru(a)$\,,
                    \item $\tru(a)\le\|a\|$\,,
                    \item $\tru$ is ultra-weakly continuous.
                \end{itemize}
                \item The standard trace
                \begin{align*}
                    \trC:\NG&\longrightarrow \C \,:\,
                    a\mapsto \langle a\cdot \delta_e , \delta_e \rangle\,.
                \end{align*}
                \item And the delocalized traces $\trg$\,, for $g\in \mathcal{C}(G)$\,,
                given by:
                \begin{align*}
                    \trg:\NG&\longrightarrow \C \,:\,
                    a\mapsto \sum_{h\in\langle g\rangle}\langle a\cdot \delta_e , \delta_h \rangle\,.
                \end{align*}
                \end{enumerate}
             \end{prop}
                These traces can be extended to $M_d(\NG)$ by taking $\tr_{\NG}^{(\cdot)}:=\tr^{(\cdot)}_{\NG}\otimes\tr_{M_d(\C)}$\,, with
                $\tr_{M_d(\C)}$ the non-normalized trace on $M_d(\C)$\,. We keep the same notation for the extended traces.
                Moreover the universal trace and the standard trace are positive,
                hence they induce the following dimension functions on $\NG$-submodules $V$ of $\ell^2(G)^d$\,:
                \begin{align*}
                    &\dim^u_{\NG}(V):=\tru\bigl(P)\in \Ze(\NG)\\
                    &\dim^\C_{\NG}(V):=\tru\bigl(P)\in \C\,.
                \end{align*}
                Here $P=P^2$ denotes the projection matrix with $\im(P)=V$\,.
                We now define the spectral density functions of operators. Roughly speaking it measures the growth of
                the spectral projections.
            \begin{defi}
                Let $G$ be a discrete group and $A\in M_d(\NG)$ a positive operator. Define
                \begin{itemize}
                     \item the spectral density function
                         \begin{align*}
                            F_A:[0,\infty)\longrightarrow [0,\infty)\,:\,
                            \varepsilon\mapsto \tr^\C_{\NG}\bigl(\chi_{[0,\varepsilon]}(A)\bigr) \,,
                        \end{align*}
                     \item and the center valued spectral density function as
                         \begin{align*}
                            F^u_A:[0,\infty)\longrightarrow\mathcal{Z}(\NG)\,:\,
                            \varepsilon\mapsto \tr^u_{\NG}\bigl(\chi_{[0,\varepsilon]}(A)\bigr) \,.
                         \end{align*}
                \end{itemize}
                where $\chi_{[0,\varepsilon]}$ denotes the characteristic function of the interval $[0,\varepsilon]$\,.
            \end{defi}

            Using this notation we have $F_A(0)=\dim^\C_{\NG}(\ker(A))$ and $F^u_A(0)=\dim^u_{\NG}(\ker(A))\,.$

            \begin{defi}
                Given $A\in\NG^+$\,, the spectral density function $F_A$ induces the Fuglede-Kadison determinant, which is defined as
                \begin{align*}
                    \lndet(A)&:=\int_{0^+}^{\infty} \ln(\lambda)\textnormal{d}F_A(\lambda)\,.
                \end{align*}
            \end{defi}

            \begin{defi}
                Let $J$ be an index set. For $A:=(a_{i,j})_{i,j\in J}$ with $a_{i,j}\in\C$\,,
                define
                \begin{align*}
                    S(A):=\sup_{i\in J}|\textnormal{supp}(z_i)|\,,
                \end{align*}
                where $z_i$ is the vector $z_i:=(a_{i,j})_{j\in J}$ and $\textnormal{supp}(z_i):=|\{j\in J\mid a_{i,j}\ne 0\}|$\,.

                Now let $|A|_\infty:=\sup_{i,j}|a_{i,j}|$ and $A^*:=(\overline{a}_{j,i})_{i,j\in J}$\,.
                Define
                \begin{displaymath}
                    \kappa(A): = \left \{
                    \begin{array}{ll} \sqrt{S(A)S(A^*)}\cdot|A|_\infty & \text{if }S(A)+S(A^*)+|A|_\infty< \infty\\
                                      \infty & \textrm{else}
                    \end{array} \right.
                \end{displaymath}
                Elements of $M_d(\C G)$ are identified with degenerated matrices, indexed by $J\times J$ where $J:=\{1,\dots,d\}\times G$\,.
                For more details we refer to \cite{Schick+}\,.
            \end{defi}

            \begin{defi}\label{defi:bound}
                Let $G$ be a  discrete group and take $A\in M_d(o(\overline{\Q})G)$ positive (where $o(\overline{\Q})$ denotes the algebraic integers),
                choose a finite Galois extension $L\subset\C$ of $\Q$\,, such that $A\in M_d(LG)$\,. Let
                $\sigma_1,\dots,\sigma_r:L\rightarrow \C$ be the different embeddings of $L$ in $\C$ with $\sigma_1$ the natural                    inclusion $L\subset\C$\,. If
                 \begin{align}\label{def:bounded}
                        \lndet(A)\ge -d\sum_{k=2}^r\ln\bigl(\kappa(\sigma_k(A))\bigl)\,,
                 \end{align}
                we say $A$ has the bounded determinant property. A discrete group $G$ is said to have the bounded determinant property,
                if all $A\in M_d(\overline{\Q}G)$ satisfy property \eqref{def:bounded}.
            \end{defi}

            \begin{lem}\label{lemma:bound}
                Given $A\in M_d(\C G)$ and let $A[i]$ be as described in \ref{sit}\,, then there exists an $i_0\in I$ such that
                for all $i\ge i_0$ we have
                \begin{align}
                    \|A|\   &\le \kappa(A)\le\infty \text{\hspace{5mm} and }\\
                    \|A[i]\|&\le \kappa(A)\,.
                \end{align}
                \begin{proof}
                    This is proven in \cite{Schick+} Lemmas 3.31\,,\,3.22\,,\,3.28\,.
                \end{proof}
            \end{lem}

            \begin{defi}\label{defi:GAE}
                Let $U<G$ be a subgroup of $G$\,. We call $G/U$ an amenable homogenous space,
                and $G$ an extension of $U$ with amenable quotient, if we have a $G$-invariant metric $d:G/U\times G/U\rightarrow \N$
                such that sets of finite diameter are finite and such that for all $K>0$ and $\varepsilon>0$ there exists some
                finite subset $\emptyset\ne X \subset G/U$ with
                \begin{align*}
                   |N_K(X)|:=|\{x\in G/U\,;\,d(x,X)\le K \text{ and } d(x,G/U-X)\le K\}|\le \varepsilon |K|\,.
                \end{align*}
                (A special case for this occurs when $U\lhd G$ is a normal subgroup and $G/U$ is an amenable group.)
            \end{defi}
            \begin{lem}\label{lemma:neig}
                A nested sequence of finite subsets $X_1\subset X_2 \subset \dots \subset G/U$ is called
                F{\o}lner exhaustion of $G/U$ if
                \,$\bigcup X_i= G/U$ and for all $K>0$ and $\varepsilon>0$ there exists an $i_0\in\N$ such that
                for all $i\ge i_0$ we have $$N_K(X_i)\le \varepsilon |X_i|\,.$$
                Every amenable homogenous space admits such an exhaustion.
                \begin{proof}
                    Compare for example Lemma 4.2 in \cite{Schick}\,.
                \end{proof}
            \end{lem}

            \section{Main Result}

               \begin{sit}
                Let $G$ be a discrete group that can be
                constructed out of groups satisfying the bounded determinant property,
                in one of the following ways:
                \begin{itemize}
                    \item $U<G$ with $\mathcal{C}(G)\subset\mathcal{C}(U)$
                    and $G/U$ admits a $G$-invariant metric making it an amenable homogenous space.
                    \item If $G$ is the direct or inverse limit of a directed system of groups $G_i$\,.
                \end{itemize}

                In $\cite{Schick+}$ the bounded determinant property is proven
                for a large class $\mathcal{G}$ of groups which is based on the above constructions.
                Most common examples with this property are amenable groups and residually finite groups.

                An other big class of groups satisfying the determinant bound property are \emph{sofic}
                groups. A brief description about \emph{sofic} groups and a proof for the determinant
                bound property is done in the next section. For more details, about \emph{sofic} groups
                we refer to \cite{Elek+Szabo} where the slightly different \emph{semi-integral-determinant}
                property is proven for sofic groups.

                We now introduce a uniform notation for the three constructions.
                Let $A\in M_d(\overline{\Q}G)$\,, where
                $\overline{\Q}$ denotes the field of algebraic numbers. The approximating matrices denoted
                by $A[i]$ will have different meanings depending on how $G$ is constructed. We have three cases.
                \begin{enumerate}\label{sit}
                    \item The group $G$ is the inverse limit of a directed system of groups $G_i$\,.
                        Define $A[i]\in M_d(\overline{\Q}G)$
                        to be the image $p_i(A)$ of $A$ under the natural map $p_i:G\rightarrow G_i$\,.
                        In this case $\tr^\C_i$\,,\,$\tr^u_{i}$ and $\tr^{\langle g\rangle}_{i}$
                        will denote $\tr^\C_{\NGI}$\,,\,$\tr^u_{\NGI}$ and $\tr^{\langle g\rangle}_{\NGI}$\,.
                    \item The group $G$ is the direct limit of a directed system of groups $G_i$\,. Denote by
                        $p_i:G_i\rightarrow G$ the corresponding maps.

                        In order to define the approximating matrices $A[i]$ we need to make some choices.
                        Write $A=(a_{k,l})$ with $a_{k,l}=\sum_{g\in G}\lambda^g_{k,l}g$\,. Then, only finitely many of the
                        $\lambda^g_{k,l}$ are non-zero. Let $V$ be the corresponding finite collection of $g \in G$\,.
                        Since $G$ is a direct limit of $G_i$\, we can find $j_0\in I$ such that $V\subset p_{j_0}(G_{j_0})$\,.
                        Choose an inverse image for each $g$ in $G_{j_0}$\,. This gives a matrix $A[j_0]\in M_d(\overline{\Q}G_{j_0})$
                        which is mapped to $A[i]:=p_{j_0i}(A[j_0])\in M_d(\overline{\Q}G_i)$ for $i>j_0$\,.
                        In this case, $\tr^\C_i$\,,\,$\tr^u_{i}$ and $\tr^{\langle g\rangle}_{i}$
                        will denote $\tr^\C_{\mathcal{N}(G_i)}$\,,\,$\tr^u_{\mathcal{N}(G_i)}$ and
                        $\tr^{\langle g\rangle}_{\mathcal{N}(G_i)}$\,.
                        Keep in mind that the values of the traces can depend on the choices made to define $A[i]$\,.
                    \item The group $G$ is an amenable extension of $U$ with F{\o}lner exhaustion
                        $X_1\subset X_2\subset \dots \subset G/U $\,.
                        Let $P_i=p_i\otimes \id_d$ with
                        $p_i:\ell^2(G)\rightarrow\ell^2(G)$ the projection on the closed subspace generated by the inverse image of
                        $X_i$ in $G$\,. The image of $P_i$ is isomorphic to $\ell^2(U)^{|X_i|d}$ as $\NU$-module.
                        We define $A[i]:=P_iAP_i$ considered as an operator on the image of $P_i$\,.

                        With this definition, $A[i]$ is no longer an element of $M_d(\NG)$ but
                        can be seen as an element in $M_{d|X_i|}(\NU)$\,. In this case,
                        $\tr^\C_i$\,,\,$\tr^u_{i}$ and $\tr^{\langle g\rangle}_{i}$ denote the following
                        \begin{align*}
                            \tr^{(\cdot)}_i(A[i])&:=\frac{1}{|X_i|}\tr^{(\cdot)}_{M_{d|X_i|}(\NU)}(A[i])\,.
                        \end{align*}
                     \end{enumerate}
                     Throughout the rest of the paper, $G_i$ will denote the obvious groups in the limit cases $(1)$ and $(2)$.
                     In the amenable case we take $G_i=U$ constantly. We use $\tr^\C_i$\,,\,$\tr^u_{i}$
                     to define $F_{A[i]}$ and $F^u_{A[i]}$\,.
            \end{sit}
            Betti-numbers are given as the dimension of the kernel of the Laplacian $\Delta_p$. Since the value of the spectral density
            functions at zero is exactly the dimension of the kernel, we can state our approximation theorem as follows.

            \begin{theorem}\label{theorem:main}
               Let $A\in M_d(\overline{\Q}G)$ and $g\in \mathcal{C}(G)$\,.
               Then, for any $\varepsilon>0$ and any choice of matrices $A[i]$\,, there exists an $i_0\in I$ such that
               for all $i\ge i_0$\,:
               \begin{align*}
                   |\langle F^u_A(0)\cdot \delta_e,\delta_g\rangle-\langle F_{A[i]}^u(0)\cdot \delta_{[e]_i},\delta_{[g]_i}\rangle|<\varepsilon\,.
               \end{align*}

                Where we denote by $\delta_{[g]_i}$ the unit vector corresponding to

                \begin{itemize}
                    \item the group element $p_i(g)\in G_i$\,, in the inverse limit case
                        $(1)$ of \eqref{sit}\,,
                    \item a chosen preimage of $g\in G_i$\,, according to the choices made to define $A[i]$ in the direct limit case $(2)$ of \eqref{sit}\,,
                    \item $g\in\mathcal{C}(G)$ in the amenable case (3) of \eqref{sit}. Without the assumption that $\mathcal{C}(G)\subset\mathcal{C}(U)$
                       approximation is still possible but then only for $g\in\mathcal{C}(U)$\,.
                \end{itemize}

            \end{theorem}

            \begin{rema}
                 The original approximation theorem (Theorem 3.12 in \cite{Schick+}) is contained in the above result if we set $g=e$\,.
            \end{rema}

            \begin{ex}
                As a direct consequence, one can use the center-valued approximation theorem
                to show the vanishing of $\beta^u$ for a closed manifold $X$ with fundamental group $\pi_1$
                in certain cases. One has
                \begin{enumerate}
                    \item $\beta^u_0(\widetilde{X})=\beta^{2}_0(\widetilde{X}) e$\,, for residually finite $\pi_1$ and
                    \item $\beta^u_p(\widetilde{X})=\beta^{2}_p(\widetilde{X}) e$,
                        for all $p\in\N$\,, if $\pi_1$ is free abelian\,.
                \end{enumerate}
                This follows directly using \cite{Lott} (example 8 and proposition 2)\,.
            \end{ex}

            \section{Bounded Determinant for Sofic Groups}

                In this section we describe the method of G. Elek and E. Sab\'{o} in \cite{Elek+Szabo} to show that
                sofic groups have the \emph{semi-integral-determinant property} and show how we can use this to prove that
                they also have the determinant bound property. We use a general method that can be used
                to show that the semi-integral determinant property implies determinant bound property if we
                have approximations with matrices over finite groups.

                \begin{defi}
                    A group $G$ has the semi-integral-determinant property if for any matrix $A\in M_d(\Z G)$ we
                    have $$\lndet(A)\ge 0\,.$$
                \end{defi}

                \begin{defi}
                    Let $G$ be a finitely generated group and $S\subset G$
                    be a finite set of generators. Then the group $G$ is called sofic, if there is a sequence of finite
                    directed graphs $\{V_n,E_n\}_{n\ge 1}$ edge-labeled by $S$ and subsets $V_0\subset V_n$ with the following
                    property:

                    For any $\delta>0$ and $r\in \N$\,, there is an integer $n_{r,\delta}$ such that if $m\ge n_{r,\delta}>0$
                    and $B_{(G,S)}(r)$ denotes the r-ball in the Cayley-graph\,,
                    then
                    \begin{itemize}
                        \item For each $v\in V_m^0$\,, there is a map $\psi:B_{(G,S)}(r)\rightarrow V_m$\,, which is an isomorphism
                        (of labeled graphs) between $B_{(G,S)}(r)$ and the r-ball in $V_m$ around $v$\,,
                        \item $|V_m^0|\ge (1-\delta)|V_m|$
                    \end{itemize}
                \end{defi}

                \begin{rema}
                    This definition for sofic groups is equivalent to the more common description using maps $\psi_n:G\rightarrow S_n$
                    and looking at the fixed-point-sets.
                \end{rema}

                \begin{theorem}\label{theorem:sofic}
                    Sofic groups have the determinant bound property (Def. \ref{defi:bound})\,.
                \end{theorem}

                Let $G$ be sofic and $A=(a_{i,j})_{1\le i,j le d} \in M_d(o(\overline{\Q}) G)$ be a positive operator. Consider the
                operator kernel of $A$\,, that is the function $K_A:G\times G\rightarrow M_d(o(\overline{\Q}))$ such that
                for $f:G\rightarrow \ell^2(G)^d$ we have
                \begin{align*}
                    Af(x)=\sum_{y\in G}K_A(x,y)f(y)\,.
                \end{align*}
                This just means $K_A(x,y)=A_{g}$ if $x=gy$ and
                $A=\sum_{g\in G} A_gg\,, A_g\in (a_{i,j}^g)_{1\le i,j \le d}\in M_d{o(\overline{\Q})}$\,.
                There is a constant $\omega_A$\,, the width of $A$ such that $K_A(x,y)=0$ if $d(x,y)>\omega_A$
                in the word metric of $G$ with respect to the generating system $S$\,.

                The approximating kernel is constructed as follows. For $m > n_{(\omega_A,\frac{1}{2})}$\,,
                define $K_A^m:V_m\times V_m\rightarrow M_d(o(\overline{\Q}))$\,, let $K_A^m(x,y)=0$ if
                $y\notin V_m^0$ and $K_A^m(x,y)=K_A^m(g,e)$ if $y\in V_m^0 \,, x=\psi_y(g)$

                \begin{lem}\label{lemma:sofic}
                    Let $G$ be a sofic group, $A\in M_d(o(\overline{\Q})G)$  a positive operator.
                    Denote by $A_m$ the bounded linear transformations on $\ell^2(V_m)^d$ defined by the kernel functions
                    $K_A^m$ and denote with $\det^*(K_A^m)$ the product of the non-zero eigenvalues of $K_A^m$\,.
                    $$\lim_{m\rightarrow \infty} \frac{\ln(\det^*(A))}{|V_m|}=\lndet(A)$$
                    \begin{proof}
                        This is proven in \cite{Elek+Szabo} Lemma (6.1)\,.
                    \end{proof}
                \end{lem}

                G. Elek and E. Szab\'{o} prove the semi-integral-determinant property (Theorem 6 in \cite{Elek+Szabo}) by using that
                the product of the positive eigenvalues of the $A_m$ are integes and hence by applying the lemma the claim follows.
                Given $A\in M_d(o(\overline{\Q})G)$\,, choose a finite Galois extension $\Q\subset L \subset \C$ such that $A\in M_d(LG)$\,.
                Let $\sigma_{i=1,\dots\,n}:L\hookrightarrow \C$ be the different embeddings of $L$ in $\C$ and denote with
                $\sigma_1$ the natural inclusion. We set $\widetilde{A}:=\bigoplus_{i=1}^d\sigma_i(A)$\,.
                For $\widetilde{A}$ Lemma \ref{lemma:sofic} obviously still holds. The product of the non-zero eigenvalues of $\widetilde{A}_m$ are
                is lowest non-zero coefficient $c$ of the characteristic polynomial. Since $o(\overline{\Q})$ is a ring,
                $c\in o(\overline{\Q})$ and $c$ is stable under all $\sigma_i$\,, $c$ is in $\Q$ and also is an algebraic integer,
                hence $c\in \Z$\,.

                \begin{lem}\label{lemma:comp1}
                    If $A$ and $B$ are positive injective operators in $M_d(\C G)$ and $A\le B$ we have
                    $$\lndet(A)\le\lndet(B)$$
                    \begin{proof}
                        This is proven in \cite{Luck})\,, Lemma 3.15\,.
                    \end{proof}
                \end{lem}

                \begin{lem}\label{lemma:comp2}
                    Let $A$ be a positive operator in $M_d(\C G)$ and let $A^\perp:\ker(A)^\perp\rightarrow \overline{\im(A)}$ be
                    the weak isomorphism obtained by restricting $A$ to $\ker(A)^\perp$\,. Then
                    $$\lndet(\sqrt{(A^\perp)^*A^{\perp}})=\lndet(A)$$
                    \begin{proof}
                        This is also proven in \cite{Luck}\,, Lemma 3.15\,.
                    \end{proof}
                \end{lem}

                We have $\sqrt{(A^\perp)^*A^{\perp}}\le \|A\|\id\le \kappa(A)\id$. By applying Lemma \ref{lemma:comp1}\,,
                Lemma \ref{lemma:comp2} and Lemma \ref{lemma:sofic} we get
                \begin{align*}
                    &0\le\lndet(\widetilde{A})=d\sum_{i=1}^n\lndet{\sigma_i(A)}\\
                    \imp&-d\sum_{i=2}^n\ln{\kappa(\sigma_i(A))}\le-\sum_{i=2}^n\lndet{\sigma_i(A)}\le \lndet(A)\,.
                \end{align*}
                This proves Theorem \ref{theorem:sofic}\,.

           \section{Some Key Lemmas}
            The Fourier coefficients of $F^u_A(0)$ are given by
             \begin{displaymath}
                     \langle F^u_A(0)\cdot \delta_e, \delta_g\rangle = \left\{
                     \begin{array}{ll}
                       \frac{1}{|\langle g\rangle|}
                            \tr^{\langle g\rangle}_{\NG}(\pr_{|\ker(A)})
                            & \text{if $g\in\mathcal{C}(G)$}\\
                        0 & \textrm{otherwise.}
                    \end{array} \right.
            \end{displaymath}
            This can be easily seen using Dixmier's approximation theorem (see e.g.
            \cite{Kadison}). In the rest of the paper $g$ is always taken in $\mathcal{C}(G)$\,.

            The proof of the $\C$-valued approximation theorem in \cite{Schick+} is based on the following three major facts.
            \begin{enumerate}
                \item $\|A\|$ and $\|A[i]\|$ have an upper bound,
                \item $\tr_{\NG}^\C$ is positive,
                \item the \emph{Fuglede-Kadison determinant} $\lndet(A)$ has a lower bound.
            \end{enumerate}
            For the center-valued approximation theorem that we prove in this paper, fact $(1)$ is obviously still valid.
            The facts $(2)$ and $(3)$ of course do not apply to our situation, since they involve the $\C$-valued trace
            $\tr^\C_{\NG}$\,, but
            the main ideas of L\"uck's method work in general for any positive functional if in addition the Fuglede-Kadison determinant
            derived from it has a lower bound for $A$ and all approximating $A[i]$\,.
            In Definition \ref{defi:help}, we define traces which are derived from delocalized traces and are positive.
            Using these traces we also define deviated Fuglede-Kadison determinants and prove the existence of a lower bound.
            Using our method, it would also be possible to directly approximate the Fourier coefficients of the projections on the homology.
            These coefficients depend on the choice of the basis, hence we do not see any application for this general approximation and
            restrict to functionals derived from delocalized traces.

            A key ingredient of our method is the following simple lemma.
            \begin{lem}
                If $a\in\NG$ is a positive element then, for all $g\in G$ we have $$\langle a\cdot \delta_e,\delta_e\rangle\ge |\langle a\cdot \delta_g,\delta_e\rangle |$$
                \begin{proof}\label{lemma:positive}
                    $a=b^*b$ then, using Cauchy-Schwarz inequality we get
    				\begin{align*}
    						\langle a\cdot \delta_e,\delta_e\rangle
    						=\|b\cdot \delta_e\|\cdot\|b\cdot \delta_e\|
    					    =\|b\cdot \delta_e\|\cdot\|b\cdot \delta_g\|
    						\ge|\langle b\cdot \delta_e,b\cdot \delta_g \rangle|
    						=|\langle a\cdot \delta_g,\delta_e\rangle |	
    				\end{align*}
    			 \end{proof}
            \end{lem}

            \begin{defi}\label{defi:help}
                Take $A\in M_d(\NG)$ and $g\in\mathcal{C}(G)-\{e\}$\,, define
                \begin{align}
                	\Tr_{\NG}^{\langle g\rangle,\re}(A):=\tr^\C_{\NG}(A)+\frac{1}{2|\langle
                    g\rangle|}\Big(\tr_{\NG}^{\langle g\rangle}(A)+\tr_{\NG}^{\langle g^{-1}\rangle}(A)\Big)\,,\\
    				\Tr_{\NG}^{\langle g\rangle,\im}(A):=\tr^\C_{\NG}(A)+\frac{1}{2i|\langle g\rangle|}\Big(\tr_{\NG}^{\langle
                    g\rangle}(A)-\tr_{\NG}^{\langle g^{-1}\rangle}(A)\Big)\,.
                \end{align}
            \end{defi}
            It follows from Lemma \ref{lemma:positive} that both traces are positive. The next lemma shows that
            for a selfadjoint $A\in M_d(\NG)$ we have
            \begin{align*}
               \langle F_A^u(0)\cdot \delta_e,\delta_g\rangle&=\frac{1}{|\langle g\rangle|}\tr_{\NG}^{\langle g\rangle}(A)\\
                &=\Tr_{\NG}^{\langle g\rangle,\re}(A)+i\Tr_{\NG}^{\langle g\rangle,\im}(A)-\tr^\C_{\NG}(A)-i\tr^\C_{\NG}(A)\,.
            \end{align*}
            This is one of the main tricks in our paper. We prove the approximation theorem for
            $\Tr_{\NG}^{\langle g\rangle,\re}$ and
            $\Tr_{\NG}^{\langle g\rangle,\im}$\,. Then we finally prove Theorem \ref{theorem:main} by
            applying this approximation and the classical approximation theorem (Theorem 3.12 in \cite{Schick+}) to the above equation.

    	    \begin{lem}\label{lemma:reel}\
    				For all $g\in\mathcal{C}(G)$ and selfadjoint $A\in M_d(\NG)$\,, the traces
                    $\Tr_{\NG}^{\langle g \rangle,\re}(A)$ and $\Tr_{\NG}^{\langle g \rangle,\im}(A)$
    				are given by the following real numbers
    				\begin{align*}
    					\Tr_{\NG}^{\langle g\rangle,\re}(A)&=\tr^\C_{\NG}(A)
                        +\re\Big(\frac{1}{|\langle g\rangle|}\tr^{\langle g\rangle}_{\NG}(A)\Big)\,,\\
    					\Tr_{\NG}^{\langle g\rangle,\im}(A)&=\tr^\C_{\NG}(A)
                        +\im\Big(\frac{1}{|\langle g\rangle|}\tr^{\langle g\rangle}_{\NG}(A)\Big)\,.	
    				\end{align*}

    				\begin{proof}
                        Since the trace on $M_d(\NG)$ is just a summation of traces on $\NG$ it is sufficient to treat the case $d=1$\,.
    					Write $A=\sum_{h\in G}\lambda_h h\in \NG$\,. We have $\langle g\rangle^{-1}=\langle g^{-1}\rangle$
                        and selfadjointness of $A$ yields $\lambda_h=\overline{\lambda_{h^{-1}}}$\,. Hence
    					\begin{align*}
    						&\tr_{\NG}^{\langle g\rangle}(A)=\overline{\tr_{\NG}^{\langle g^{-1} \rangle}(A)}\,.
    					\end{align*}

    				\end{proof}

    	   \end{lem}

        \section{Lower Bound for Determinants}

            \begin{defi}
                Take a positive operator $A\in M_d(\NG)$ and denote by
                $\{E^A_\lambda:=\chi_{[0,\lambda]}(A)\mid\lambda\in \R^+_0\}$  the spectral
                family of $A$\,. Define then the following spectral density functions:
                \begin{align*}
                    &F_A(\lambda):=\tr^\C_{\NG}(E^A_\lambda)\,,\\
                    &F^{\langle g\rangle,\re}_A(\lambda):=\Tr_{\NG}^{\langle g\rangle,\re}(E^A_\lambda)\,,\\
                    &F^{\langle g\rangle,\im}_A(\lambda):=\Tr_{\NG}^{\langle g\rangle,\im}(E^A_\lambda)\,.
                \end{align*}
            \end{defi}

            For positive $A\in M_d(\NG)$\,, the spectral density functions
            $F_A\,,\,F^{\langle g\rangle,\re}_A$ and $F^{\langle g\rangle,\im}_A$ are monotone increasing and induce
            Riemann-Stieltjes measures
            $\textnormal{d}F_A(\lambda)\,,\,\textnormal{d}F_A^{\langle g\rangle,\re}(\lambda)$
            and $\textnormal{d}F_A^{\langle g\rangle,\im}(\lambda)$\,,
            allowing us to define the following (deviations of the) Fuglede-Kadison determinant.

            \begin{defi}
                Take $A\in M_d(\NG)$ positive and define
                \begin{align*}
                    \lndet(A)&:=\int_{0^+}^{\infty} \ln(\lambda)\textnormal{d}F_A(\lambda)\,,\\
                    \lndet^{\langle g\rangle,\re}(A)&:=\int_{0^+}^{\infty} \ln(\lambda)\textnormal{d}F_A^{\langle g\rangle,\re}(\lambda)\,,\\
                    \lndet^{\langle g\rangle,\im}(A)&:=\int_{0^+}^{\infty} \ln(\lambda)\textnormal{d}F_A^{\langle g\rangle,\im}(\lambda)\,.
                \end{align*}
            \end{defi}

            In order to prove Theorem $\ref{theorem:main}$ we need a lower bound for the deviated Fuglede-Kadison determinants
            $\lndet^{\langle g\rangle,\re}(A)$ and $\lndet^{\langle g\rangle,\re}(A)$\,. We obtain it using the fact that the perturbation
            caused by the delocalized trace is controlled by the standard trace.

            \begin{lem}\label{lemma:lbound}
                          Let $G$ be a group that satisfies the determinant bound property and is constructed as described in \ref{sit}\,.
                          Take $A\in M_d(o(\overline{\Q})G)$ positive (where $o(\overline{\Q})$ denotes the algebraic integers),
                          choose a finite Galois extension $L\subset\C$ of $\Q$\,, such that $A\in M_d(LG)$\,. Let
                          $\sigma_1,\dots,\sigma_r:L\rightarrow \C$ be the different embeddings of $L$ in $\C$ with $\sigma_1$ the natural
                          inclusion $\sigma_1:L\subset\C$\,.
                          Then
                          \begin{align*}
                            \lndet^{\langle g\rangle,\re}(A)\ge -2d\big|\sum_{k=2}^r\ln\bigl(\kappa(\sigma_k(A))\bigl)\big|\,,\\
                            \lndet^{\langle g\rangle,\im}(A)\ge -2d\big|\sum_{k=2}^r\ln\bigl(\kappa(\sigma_k(A))\bigl)\big|\,.
                         \end{align*}

                 \begin{proof}
                    We prove the lemma only for $\lndet^{\langle g\rangle,\re}$\,, the case $\lndet^{\langle g\rangle,\im}$ being identical. Using Lemmas \ref{lemma:positive} and \ref{lemma:reel} we get
                    \begin{align}\label{eq:lbound}
                        \tr_{\NG}^\C(A)\ge\frac{1}{|\langle g\rangle|}|\re\bigr(\tr_{\NG}^{\langle g\rangle}(A)\bigl)|\,.
                    \end{align}
                    Define the function
                    $f_A^{\langle g\rangle,\re}(\lambda):=\frac{1}{|\langle g\rangle|}\re\bigl(\tr_{\NG}^{\langle g\rangle}(E_\lambda^A))$\,.
                    For $a\le b\in\R_0^+$\,, inequality \eqref{eq:lbound} yields
                    \begin{align*}
                        F_A(b)-F_A(a)\ge |f_A^{\langle g\rangle,\re}(b)-f_A^{\langle g\rangle,\re}(a)|\,.
                    \end{align*}
                    The Riemann-Stieltjes measure induced by $F_A(\lambda)$ dominates in absolute values
                    the (possibly signed) measure induced by $f_A(\lambda)$\,. Hence, we have
                    \begin{align*}
                        \lndet^{\langle g\rangle,\re}(A)&=\int_{0^+}^{\infty} \ln(\lambda)\textnormal{d}F_A^{\langle g\rangle,\re}(\lambda)\\
                            &=\int_{0^+}^{\infty} \ln(\lambda)\textnormal{d}F_A(\lambda)+
                                \int_{0^+}^{\infty} \ln(\lambda)\textnormal{d}f_A^{\langle g\rangle,\re}(\lambda)\\
                            &\ge-\big|\int_{0^+}^{\infty} \ln(\lambda)\textnormal{d}F_A(\lambda)\big|
                               -\big|\int_{0^+}^{\infty} \ln(\lambda)\textnormal{d}f_A^{\langle g\rangle,\re}(\lambda)\big|\\
                            &\ge-2\big|\int_{0^+}^{\infty} \ln(\lambda)\textnormal{d}F_A(\lambda)\big|\\
                            &\ge -2d\big|\sum_{k=2}^r\ln\bigl(\kappa(\sigma_k(A))\bigl)\big|
                    \end{align*}
                 \end{proof}
            \end{lem}

    \section{Convergence of the Trace}

        In this section we basically use the ideas from
        \cite{Schick+} and \cite{Schick} to prove the following equalities, for all $G$ constructed as described in Situation \ref{sit}\,, $g\in\mathcal{C}(G)$\,,
        $A\in M_d(\C G)$ and every polynomial $p\in\C[x]$\,:
        \begin{align}\label{eq:limit1}
            \lim_{i\rightarrow \infty}\Tr^{\langle g\rangle,\re}_i(p(A[i]))&=\Tr_{\NG}^{\langle g\rangle,\re}(p(A))\,,\\
            \lim_{i\rightarrow \infty}\Tr^{\langle g\rangle,\im}_i(p(A[i]))&=\Tr_{\NG}^{\langle g\rangle,\im}(p(A))\,.\label{eq:limit2}
        \end{align}
        The traces $\Tr_i$ depend on the construction of $G$\,.
        We deal first with the limit cases $(1)$ and $(2)$ of \ref{sit}\,.
        \begin{lem}
            Take
            $A\in M_d(\C G)$\,,\,$p\in\C[x]$ and $g\in\mathcal{C}(G)$\,.
            If $G$ is the direct or inverse limit of groups $(G_i)_{i\in I}$ then there is an $i_0\in I$ such that
            for all $i\ge i_0$\,:
            \begin{align*}
                \Tr^{\langle g\rangle,\re}_i(p(A[i]))&=\Tr_{\NG}^{\langle g\rangle,\re}(p(A))\,,\\
                \Tr^{\langle g\rangle,\im}_i(p(A[i]))&=\Tr_{\NG}^{\langle g\rangle,\im}(p(A))\,.
            \end{align*}
            \begin{proof}
                The proof follows directly from the fact that the \emph{support}
                \begin{align*}
                    \textnormal{supp}(p(A)):=\Big\{\lambda_g^{k,l}\ne 0\mid 1\le k,l \le d\,,\,(p(A))_{k,l}=\sum_{g\in G} \lambda_g^{k,l}g\Big\}
                \end{align*}
                of $p(A)\in M_d(\C G)$ is finite.
                Since $G$ is an inverse or direct limit, choosing $i_0$ big enough, we have, for all $i\ge i_0$:
                \begin{align*}
                    \textnormal{supp}(p(A[i]))=\textnormal{supp}(p(A))\,.
                \end{align*}
                As a consequence, the traces coincide.
            \end{proof}
        \end{lem}

        To prove  \eqref{eq:limit1} and \eqref{eq:limit2} in
        the amenable case (3) of \ref{sit}, we adapt ideas from \cite{Schick}\, (Lemma 4.6) to our situation.
        \begin{lem}
            Let $G$ be an amenable extension of $U$ with F{\o}lner exhaustion $X_1\subset X_2 \subset\dots\subset G/U$\,.
            Then, for all $g\in\mathcal{C}(U)$\,, $A\in M_d(\C G)$ and every polynomial $p\in\C[x]$ we have
            \begin{align*}
                \lim_{i\rightarrow \infty}\Tr^{\langle g\rangle,\re}_i(p(A[i]))&=\Tr_{\NG}^{\langle g\rangle,\re}(p(A))\,,\\
                \lim_{i\rightarrow \infty}\Tr^{\langle g\rangle,\im}_i(p(A[i]))&=\Tr_{\NG}^{\langle g\rangle,\im}(p(A))\,.
            \end{align*}

            \begin{proof}
                Again we only treat the case $\Tr_{\NG}^{\langle g\rangle,\re}$ and assume $d=1$\,, since the
                general case follows by summing up the traces. Let $A\in\NG$ and denote
                $A[i]:=P_i A P_i^*$\,, as described in \ref{sit}\,. By linearity of the trace, it
                also suffices to treat the case where $p$ is a monomial. Pull back the metric on $G/U$
                in order to get a semi-metric on $G$\,. Denote the inverse image of $X_i$ by $X_i'$\,.
                For $g\in X_i'$ and $h\in U$we have $P_i(h\cdot\delta_g)=h\cdot\delta_g$\,. Selfadjointness of $P_i$
                 implies for $g\in X_i'$ $h\in U$\,, that $\langle(P_iAP_i)^n \delta_g,h\cdot\delta_g\rangle=\langle AP_i AP_i\dots P_iA\delta_g,h\cdot\delta_g\rangle$ and we have the following telescope sum:
                  \begin{align}\label{eq:tele}
                    AP_iA \cdots P_i
                    A=A^n-A(1-P_i)A^{n-1}
                    \cdots-AP_i\cdots A(1-P_i)A\,.
                  \end{align}
                  We now compute for $s\in \mathcal{C}(U)$\,,
                  \begin{align*}
                    \Big|\Tr^{\langle s\rangle,\re}_{\NG}(A^n)-\Tr^{\langle s\rangle,\re}_{i}(A[i]^n)\Big|
                    =&\Big|\langle A_n\delta_e,\delta_e\rangle
                        +\frac{1}{2|\langle s\rangle|}\sum_{h\in\langle s\rangle\cup
                        \langle s^{-1}\rangle}\langle A^n\delta_e,h\cdot\delta_e\rangle\\
                        &-\frac{1}{|X_i|}\sum_{[g]\in X_i}\Big(\langle A^n\delta_g,\delta_g\rangle
                            -\frac{1}{2|\langle s\rangle|}\sum_{h\in\langle
                            s\rangle\cup\langle s^{-1}\rangle}\langle A[i]^n\delta_g,h\cdot\delta_g\rangle\Big)\Big|\\
                    \le&\frac{1}{|X_i|}\sum_{[g]\in X_i}
                            \Big|\langle A_n\delta_g,\delta_g\rangle
                        +\frac{1}{2|\langle s\rangle|}\sum_{h\in\langle s\rangle\cup
                        \langle s^{-1}\rangle}\langle A^n\delta_g,h\cdot\delta_g\rangle\\
                        &-\langle A^n\delta_g,\delta_g\rangle
                            -\frac{1}{2|\langle s\rangle|}\sum_{h\in\langle
                            s\rangle\cup\langle s^{-1}\rangle}\langle A[i]^n\delta_g,h\cdot\delta_g\rangle\Big|\\
               \end{align*}
               \begin{align*}
                    =&\frac{1}{|X_i|}\sum_{[g]\in X_i}
                            \Big|\Bigl(\langle A^n\delta_g,\delta_g\rangle
                            -\langle A^n\delta_g,\delta_g\rangle\Bigr)
                        +\frac{1}{2|\langle s\rangle|}
                            \sum_{h\in\langle s\rangle\cup\langle s^{-1}\rangle}\Bigl(\langle A^n \delta_g,h\cdot \delta_g\rangle
                            -\langle A[i]^n \delta_g,h\cdot\delta_g\rangle\Bigr)\Big|\,.
                  \end{align*}
                  Using \eqref{eq:tele} and applying Cauchy-Schwartz inequality, we get
                  \begin{align*}
                    \Big|\Tr^{\langle g\rangle,\re}_{\NG}(A^n)-\Tr^{\langle g\rangle,\re}_{i}(A[i]^n)\Big|
                    \le&\frac{1}{|X_i|}\sum_{j=1}^{n-1}\sum_{[g]\in X_i}\Big|
                            \langle (1-P_i)A^j\delta_g,(A^*P_i)^{n-j}\delta_g\rangle\\
                        &+\frac{1}{2|\langle s\rangle|} \sum_{h\in\langle s\rangle\cup
                            \langle s^{-1}\rangle}\langle (1-P_i)A^j\delta_g,(A^*P_i)^{n-j} h\cdot \delta_g\rangle\Big|\\
                    \le&\frac{1}{|X_i|}\sum_{j=1}^{n-1}\sum_{[g]\in X_i/U}\Bigl(
                            \|(1-P_i)A^j\delta_g\|\cdot\|A^*\|^{n-j}\\
                        &+\frac{1}{2|\langle s\rangle|}
                            \sum_{h\in\langle s\rangle\cup\langle s^{-1}\rangle}\|(1-P_i)A^j\delta_g\|\cdot\|A^*\|^{n-j}\Bigr)\\
                    \le&\frac{2}{|X_i|}\sum_{j=1}^{n-1}\sum_{[g]\in X_i}
                            \|(1-P_i)A^j\delta_g)\|\cdot\|A^*\|^{n-j}\,.
                  \end{align*}
                  Define for $i\in \N$ $$T_i:=\Big\{g\in G \,\,|\,\, \lambda_{[i],g}^{k,l}\ne 0 \text{ where }
                        \,(A[i])_{k,l}:=\sum_{g\in G}\lambda_{[i],g}^{k,l}g \text{ and } 1\le k,l\le d\Big\}\,.$$
                  Then the set $T:=\bigcup\limits_{i=1}^{\infty}T_i$ is a finite subset of $G$\,.
                  Hence if we take $R\in\N$ big enough and let $B_R(g)$ be the ball with radius $R$ around $g$\,, we have
                  \begin{align*}
                        (1-P_{B_R(g)})A^j\delta_g=0\,.
                  \end{align*}
                  The integer $R$ is independent from $g$\,, since the semi-metric is $G$ invariant.
                  Now if $B_R(g)\subset X_i'$\,, which means $[g]\in X_i-N_R(X_i)$ (see Definition \ref{defi:GAE})\,,
                  we have  $\im(P_{B_R})\subset\im(P_i)$ and hence
                  \begin{align*}
                        (1-P_i)A^j\delta_g=0\,.
                  \end{align*}
                  Now we have
                  \begin{align*}
                    \Big|\Tr^{\langle s\rangle,\re}_{\NG}(A^n)-\Tr^{\langle s\rangle,\re}_{i}(A[i]^n)\Big|
                    &\le \frac{2}{|X_i|}\sum_{j=1}^{n-1}\sum_{[g]\in X_i}
                        \|(1-P_i)A^j\delta_g\|\cdot\|A^*\|^{n-j}\\
                    &=\frac{2}{|X_i|}\sum_{j=1}^{n-1}\sum_{[g]\in N_R(X_i)}
                        \|(1-P_i)A^j\delta_g\|\cdot\|A^*\|^{n-j}\\
                    &\le\frac{|N_R(X_i)|}{|X_i|}
                    2\sum_{j=1}^{n-1}\|(1-P_i)A^j\|\cdot\|A^*\|^{n-j}\\
                    &\le\frac{|N_R(X_i)|}{|X_i|}
                    \underbrace{2n\max_{j=1,\dots,n}\{\|A\|^j\cdot\|A^*\|^{n-j}\}}_{c_n}.
                  \end{align*}
                  The quantity $c_n$ is independent of $i$ and
                  Lemma \ref{lemma:neig} shows that $$\lim_{i\rightarrow\infty}\frac{|N_R(X_i)|}{|X_i|}=0\,;$$
                  hence the claim follows.
            \end{proof}

        \end{lem}
    \section{Finalization of the Proof}

            Now we are finally ready to prove our theorem. The
            main idea in this section is to use the lower bound of the Fulglede-Kadision determinant
            and is due to W. L\"uck in \cite{Luck1}\,.

            Define for the spectral density functions $F^{\langle g\rangle,\re}_A$ and $F^{\langle g\rangle,\im}_A$:
                \begin{align*}
                    \overline{F}^{(\cdot)}_A(\lambda):=\limsup\limits_{i\rightarrow\infty}(F^{(\cdot)}_{A[i]})(\lambda)\,,\\
                    \underline{F}^{(\cdot)}_A(\lambda):=\liminf\limits_{i\rightarrow\infty}(F^{(\cdot)}_{A[i]})(\lambda)\,,
                \end{align*}
                and denote their right-continuous approximations by
                \begin{align*}
                    \overline{F}^{(\cdot),+}_A(\lambda):=\lim\limits_{\varepsilon\rightarrow 0^+}(\overline{F}^{(\cdot)}_{A})(\lambda+\varepsilon)\,,\\
                    \underline{F}^{(\cdot),+}_A(\lambda):=\lim\limits_{\varepsilon\rightarrow 0^+}(\underline{F}^{(\cdot)}_{A})(\lambda+\varepsilon)\,.
                \end{align*}

            \begin{theorem}
                Let $g\in\mathcal{C}(G)$ and $A\in M_d(\overline{\Q} G)$\,. Then
                \begin{align*}
                    & F^{\langle g\rangle,\re}_A(0)=\lim_{i\rightarrow\infty}F^{\langle g\rangle,\re}_{A[i]}(0)\,,\\
                    & F^{\langle g\rangle,\im}_A(0)=\lim_{i\rightarrow\infty}F^{\langle g\rangle,\im}_{A[i]}(0)\,.
                \end{align*}

                \begin{proof}
                    We  only prove $F^{\langle g\rangle,\re}_A(0)=\lim_{i\rightarrow\infty}F^{\langle g\rangle,\re}_{A[i]}(0)$\,. The other case
                    can be done identically.

                    Fix $\lambda\ge 0$ and take a sequence $P_n$ of polynomials converging
                    pointwise to $\chi_{[0,\lambda]}$\,, such that for $0\le x \le \kappa(A)$\,,
                    \begin{align*}
                        \chi_{[0,\lambda]}(x)\le &P_n(x) \le \chi_{[0,\lambda+\frac{1}{n}]}(x)+\frac{1}{n}\chi_{[0,\kappa(A)]}(x)\,.
                        \intertext{Applying functional calculus preserves the inequality and
                        since for all $i\in I$ $\|A[i]\|\le\kappa{(A)}$ we get}
                        E^{A[i]}_{\lambda}\le &P_n(A[i]) \le E^{A[i]}_{\lambda+\frac{1}{n}}+\frac{1}{n}\id\,.
                        \intertext{Then we apply the positve and hence order preserving trace $\trrei$ and
                        use the fact that $\trrei(\id)\le2\tr_{\NG}^\C(\id)=2d$\,. We get
                                for all $i\in I$}
                         F^{\langle g\rangle,\re}_{A[i]}(\lambda)\le &\trrei(P_n(A[i]))\le F^{\langle
                           g\rangle,\re}_{A[i]}(\lambda+\frac{1}{n})+\frac{2d}{n}\,.
                        \intertext{Taking $\limsup$ on the left side and $\liminf$ on the right side leads to:}
                        \overline{F}^{\langle g\rangle,\re}_{A}(\lambda)\le &\trre(P_n(A))\le\underline{F}^{\langle g\rangle,\re}_A(\lambda+\frac{1}{n})+\frac{2d}{n}\,.
                    \end{align*}
                    The sequence $P_n(A)$ converges strongly in a norm bounded set, hence it converges already in the ultra-strong topology.
                    Taking $n\rightarrow\infty$ and using normality of $\trrei$ yields:
                    \begin{align}\label{eq:abschaetzung}
                        \overline{F}^{\langle g\rangle,\re}_{A}(\lambda)\le
                            &F^{\langle g\rangle,\re}_A(\lambda)\le\underline{F}^{\langle g\rangle,\re,+}_A(\lambda)\,.
                    \end{align}
                    Setting $\lambda=0$ gives us the first half of the proof:
                    \begin{align*}
                        \limsup_{i\in I}{F}^{\langle g\rangle,\re}_{A[i]}(0)\le &F^{\langle g\rangle,\re}_A(0)\,.
                    \end{align*}
                    We prove now that
                    $F^{\langle g\rangle,\re}_A(0)\le\liminf_{i\in I}{F}^{\langle g\rangle,\re}_{A[i]}(0)$\,, which finishes the proof.
                    We first pass from $I$ to a subnet $J\subset I$\,, such that $\limsup_{i\in J}F^{\langle g\rangle,\re}_A(0)
                    =\liminf_{i\in I}F^{\langle g\rangle,\re}_A(0)$\,. Equation \eqref{eq:abschaetzung} still holds and we keep our notation
                    $\overline{F}^{(\cdot)}\,,\,\underline{F}^{(\cdot)}$ but using $J$ instead of $I$\,.
                    Moreover we need the Fatou lemma and the fact that the (deviated) Fuglede-Kadison determinant is bounded.
                    For this we restrict to the case $A\in M_d(\overline{\Q}G)$
                    to the case $M_d(o(\overline{\Q})G)$\,, since Lemma \ref{lemma:lbound} only holds for
                    $A\in M_d(o(\overline{\Q})G)$\,. But every algebraic number $z$ can be written as
                    a quotient $y/k$ with $y\in o(\overline{\Q})$ and $k\in \N$\,. We work then with $s A\in M_d(o(\overline{\Q})G)$
                    instead of $A\in M_d(\overline{\Q}G)$, where $s$ is an appropriate integer.
                    Of course this does not change the kernel and we do not lose any generality.

                    Recall that $\kappa(A)\ge \|A\|,\|A[i]\|$\,. Using partial integration, we get
                    \begin{align*}
                        &\lndet^{\langle g\rangle,\re}(A)=\ln({\kappa(A)})({F}^{\langle g\rangle,\re}_{A}(\lambda)-{F}^{\langle g\rangle,\re}_{A}(0))-
                            \int_{0^+}^{\kappa(A)}\frac{{F}^{\langle g\rangle,\re}_{A}(\lambda)-{F}^{\langle
                        g\rangle,\re}_{A}(0)}{\lambda}\textnormal{d}\lambda
                    \end{align*}
                    Lemma \ref{lemma:lbound} yields a $C\in \R$\,, independent of $i\in I$\,, such that
                    $\lndet^{\langle g\rangle,\re}(A[i])\ge C$
                    and since ${F}^{\langle g\rangle,\re}_{A[i]}(\lambda)\le\trrei(\id)\le 2d$\,, it follows that
                    \begin{align}\label{eq:positive}
                        &\int_{0^+}^{\kappa(A)}\frac{{F}^{\langle g\rangle,\re}_{A[i]}(\lambda)-{F}^{\langle
                        g\rangle,\re}_{A[i]}(0)}{\lambda}\textnormal{d}\lambda
                        \le \ln({\kappa(A)})({F}^{\langle g\rangle,\re}_{A[i]}(\lambda)-{F}^{\langle g\rangle,\re}_{A[i]}(0))\le
                        2d\cdot\ln({\kappa(A)})-C
                    \end{align}

                    Moreover for $\varepsilon\ge 0$ we get
                    \begin{align*}
                        &\Big|\int_{\varepsilon}^{\kappa(A)}\frac{\underline{F}^{\langle g\rangle,\re ,+}_{A}(\lambda)-{F}^{\langle
                        g\rangle,\re}_{A}(0)}{\lambda}\textnormal{d}\lambda
                        -\int_{\varepsilon}^{\kappa(A)}\frac{\underline{F}^{\langle g\rangle,\re}_{A}(\lambda)-{F}^{\langle
                        g\rangle,\re}_{A}(0)}{\lambda}\textnormal{d}\lambda\Big|\\
                        =&\lim_{n\rightarrow\infty}\Big|\int_{\varepsilon}^{\kappa(A)}\frac{\underline{F}^{\langle
                        g\rangle,\re}_{A}(\lambda+\frac{1}{n})
                        -{F}^{\langle g\rangle,\re}_{A}(0)}{\lambda}\textnormal{d}\lambda
                        -\int_{\varepsilon}^{\kappa(A)}\frac{\underline{F}^{\langle g\rangle,\re }_{A}(\lambda)-{F}^{\langle
                        g\rangle,\re}_{A}(0)}{\lambda}\textnormal{d}\lambda\Big|\\
                        =&\lim_{n\rightarrow\infty}\Big|\int_{\varepsilon+\frac{1}{n}}^{{\kappa(A)}+\frac{1}{n}}
                        \frac{\underline{F}^{\langle g\rangle,\re}_{A}(\lambda)-{F}^{\langle g\rangle,\re}_{A}(0)}{\lambda}\textnormal{d}\lambda
                        -\int_{\varepsilon}^{\kappa(A)}\frac{\underline{F}^{\langle g\rangle,\re}_{A}(\lambda)-{F}^{\langle
                        g\rangle,\re}_{A}(0)}{\lambda}\textnormal{d}\lambda\Big|\\
                        =&\lim_{n\rightarrow\infty}\Big|\int_{{\kappa(A)}}^{{\kappa(A)}+\frac{1}{n}}
                        \frac{\underline{F}^{\langle g\rangle,\re}_{A}(\lambda)-{F}^{\langle g\rangle,\re}_{A}(0)}{\lambda}\textnormal{d}\lambda
                        -\int_{\varepsilon}^{\varepsilon+\frac{1}{n}}\frac{\underline{F}^{\langle g\rangle,\re}_{A}(\lambda)-{F}^{\langle
                        g\rangle,\re}_{A}(0)}{\lambda}\textnormal{d}\lambda\Big|\\ \le&\lim_{n\rightarrow\infty}\Bigl(\frac{\underline{F}^{\langle
                        g\rangle,\re}_{A}({\kappa(A)})-{F}^{\langle g\rangle,\re}_{A}(0)}{n\varepsilon}
                        -\frac{\underline{F}^{\langle g\rangle,\re}_{A}({\kappa(A)})-{F}^{\langle g\rangle,\re}_{A}(0)}{n{\kappa(A)}}\Bigr)=0
                    \end{align*}
                    Since this holds for every $\varepsilon>0$ we can now use equation \eqref{eq:positive} to finish the proof
                    \begin{align*}
                        \int_{0^+}^{\kappa(A)}\frac{{F}^{\langle g\rangle,\re}_{A}(\lambda)-{F}^{\langle
                            g\rangle,\re}_{A}(0)}{\lambda}\textnormal{d}\lambda
                        \le&\int_{0^+}^{\kappa(A)}\frac{\underline{F}^{\langle g\rangle,\re ,+}_{A}(\lambda)-{F}^{\langle
                            g\rangle,\re}_{A}(0)}{\lambda}\textnormal{d}\lambda\\
                        =&\int_{0^+}^{\kappa(A)}\frac{\underline{F}^{\langle g\rangle,\re}_{A}(\lambda)-{F}^{\langle
                            g\rangle,\re}_{A}(0)}{\lambda}\textnormal{d}\lambda\\
                        \underset{(*)}{\le}&\int_{0^+}^{\kappa(A)}\frac{\underline{F}^{\langle g\rangle,\re}_{A}(\lambda)-\overline{F}^{\langle
                            g\rangle,\re}_{A}(0)}{\lambda}\textnormal{d}\lambda\\
                        \le&\int_{0^+}^{\kappa(A)}\frac{\liminf_{i\in J}\Big({F}^{\langle g\rangle,\re}_{A[i]}(\lambda)-{F}^{\langle
                            g\rangle,\re}_{A[i]}(0)\Big)}{\lambda}\textnormal{d}\lambda\\
                        \le&\liminf_{i\in J}\int_{0^+}^{\kappa(A)}\frac{{F}^{\langle g\rangle,\re}_{A[i]}(\lambda)-{F}^{\langle
                            g\rangle,\re}_{A[i]}(0)}{\lambda}\textnormal{d}\lambda\\
                        \le& 2d\ln({\kappa(A)})
                    \end{align*}
                    From this it follows that $\underline{F}^{\langle g\rangle,\re, +}_{A}(0)=\overline{F}^{\langle g\rangle,\re}_{A}(0)$\,, otherwise the third integral $(*)$ would not be finite. So we have, using equation \eqref{eq:abschaetzung}
                    \begin{align*}
                        \liminf_{i\in I}F^{\langle g\rangle,\re}_{A[i]}(0)=\limsup_{i\in J}F^{\langle g\rangle,\re}_{A[i]}(0)=F^{\langle g\rangle,\re}_{A}(0)\,,
                    \end{align*}
                    hence the second part is proven.
                \end{proof}

            \end{theorem}

            \nocite{Luck}
            \nocite{Atiyah}
            \nocite{Clair}
            \nocite{Farber}
            \nocite{Linnell}
            \nocite{Luck2}
            \bibliographystyle{plain}
         \bibliography{bib}
\end{document}